\documentclass{article}

\usepackage{graphicx}
\usepackage{subfig}
\usepackage{natbib}
\usepackage{bm}
\usepackage{amsmath, amsthm, amsfonts, amsbsy, amssymb, upref} 
\usepackage{caption}

\newcommand{\be}{\begin{equation}}
\newcommand{\ee}{\end{equation}}
\newcommand{\mR}{\mathcal{R}}

\renewcommand{\vec}[1]{\bm{#1}}

\newcommand{\xs}{x_{1}}

\begin{document}

\title{A Selective Approach to Internal Inference}

\author{Samuel M. Gross\thanks{email:smgross@stanford.edu}  \and Jonathan Taylor\thanks{email:Supported by NSF Grant
    NSF grant DMS-12-08857} \and Robert
  Tibshirani\thanks{email:tibs@stanford.edu, Supported by NSF Grant
    DMS-99-71405 and National Institutes of Health Contract
    N01-HV-28*183}}

\maketitle

\begin{abstract}
A common goal in modern biostatistics is to form a \emph{biomarker
  signature} from high dimensional gene expression data that is
predictive of some outcome of interest.  After learning this biomarker
signature, an important question to answer is how well it predicts the
response compared to classical predictors.  This is challenging, because the biomarker signature is an \emph{internal predictor} -- one that has been learned using the same dataset on which we want to evaluate it's significance.  We propose a new method for approaching this problem based on the technique of selective inference.  Simulations show that our method is able to properly control the level of the test, and that in certain settings we have more power than sample splitting.
\end{abstract}

\section{Introduction}

Since the proliferation of  microarrays and other high throughput assays, there have been many studies seeking to find  biomarker signatures that are predictive of either disease status or outcomes under various treatments.  One natural and important question upon discovering these  biomarker signatures is whether they provide predictive value beyond that of standard clinical predictors such as height, weight, and age.  This is a challenging task, since the signature was discovered on the same dataset that is to be used for performing inference and this gives the it an unfair advantage.  We seek a statistical approach that makes a fair inference in this setting and also exhibits good power.

To make our setting more concrete, suppose that we have a  $n\times p_x$ matrix of expression data $X$ with $n$ samples
and $p_x$ biomarkers (features), a quantitative
outcome variable $y$ of length $n$, and an $n\times p_z$ matrix of clinical factors
$Z$.   We assume that $p_z << n << p_x$, so that $Z$ would be considered a low dimensional dataset and $X$ a high dimensional dataset.  

At the first stage, we derive a biomarker signature $\hat y$ of length $n$ from $(X,y)$.  We call these signatures {\em internal} predictors since they are derived from the present dataset.  In contrast, we call the columns
$Z$  {\em external} predictors since they are not a function of our dataset.

At the second stage of our analysis, we might analyze $\hat y$ by fitting a regression model of $y$ on $\hat y$ and $Z$ and evaluating the significance of the coefficient for $\hat y$.
But there is an obvious problem with this: $\hat y$ was derived from the same data being used to
calculate significance and hence has an unfair advantage. In particular, an overfit $\hat y$
will look unreasonably good in the regression of $y$ on $\hat y$ and $Z$.  This is even more of a concern because $X$ is high dimensional, so typically a model is available where $\hat{y} \approx y$.

Formally, we want to perform inference in a two step procedure.  Assume that we have some function $f(X_2; X_1, y_1)$ that returns the fit of $X_2$ using parameter estimates learned from $(X_1, y_1)$.  For example, in the case of Ordinary Regression: 
\[f_{OR}(X_2; X_1, y_1) = X_2 (X_1^TX_1)^{-1}X_1^Ty_1.\]  
We will sometimes refer to these functions as {\em fitters}.  In the first step, we build our internal predictor: 
\[\hat{y} \equiv f(X;X,y).\]

Since we are in a high dimensional setting ($n <<< p_x$), we focus our attention on fitters that involve some sort of variable selection.  Namely, $f(X_2; X_1, y_1)$ should depend on only some of the columns of $X_2$ (exactly which columns will depend on $(X_1, y_1)$).  We call the set of columns of $X_2$ that are selected by our fitter $\hat{E}$, so that $(X_2)_{\hat{E}}$ are those columns as a matrix.

In the second stage, we consider two different ways of evaluating our internal predictor.  The first way is to directly test the inclusions of $\hat{y}$ in a regression of $y$ on $Z$.  Although our test statistic is not generically Student-t distributed in this case, this test is similar to the typical t-test for the inclusion of a variable in regression.  The other way we consider is to test the inclusion of the predictors $X_{\hat{E}}$ in a regression of $y$ on $Z$.  This test is similar to an F-test, although again, our test statistic is not generically F distributed.  The two tests we consider are:

\begin{enumerate}
\item
\emph{Test the inclusion of $\hat{y}$:}
\begin{equation}
  H_0:\theta = 0,\;\;\;y = \theta\hat{y} + Z\gamma + \epsilon
\end{equation}
where $\gamma \in \mR^{n\times p_z}$ is a nuisance parameter, and $\epsilon$ is an iid normally distributed noise vector with variance $\sigma^2$.

\item
\emph{Test the inclusion of $X_{\hat{E}}$:}
\begin{equation}
  H_0:\beta = 0,\;\;\;y =  X_{\hat{E}}\beta + Z\gamma + \epsilon
\end{equation}
where $\gamma \in \mR^{n\times p_z}$ is a nuissance parameter, and $\epsilon$ is an iid normally distributed noise vector with variance $\sigma^2$.
\end{enumerate}

It is important to realize that these two tests are testing different things.  The test for the inclusion of $X_{\hat{E}}$ allows one to find a new direction in $X_{\hat{E}}$ that does a good job of explaining $y - Z\hat{\gamma}$.  On the other hand, the test for the inclusion of $\hat{y}$ never sees $Z$ when it builds $\hat{y}$.  The test for the inclusion of $\hat{y}$ is more similar to how people perform internal inference in the past, but the test for the inclusion of $X_{\hat{E}}$ may have value as well.

To give the reader an idea of the issue here, we present some results from a simulation where $X$ is independent of both $Z$ and $y$.  Results for one instance, as well as the type I error over 2000 instances are summarized in Table~\ref{tab:null}.  In the first stage we build a lasso model with fixed $\lambda$ to fit $\hat{y}$ \citep{atibshirani1996regression}:
\[\hat\eta = \arg\min\frac{1}{2}\|y-X\eta\|_2^2 + \lambda\|\eta\|_1, \;\;\; \hat{y} = X\hat{\eta}.\]
In the second stage, we test the inclusion of $\hat{y}$ in the regression of $y$ on $Z$.  All of the tests were run to control type-I error at .05.

\begin{table}
\centering
\caption{Null Simulation with $\alpha=0.05$}
\begin{tabular}{ | l | c | c|}
\hline
  Method & p-value & Actual type I error\\
\hline \hline
  Naive T-test & 1.69e-07 & 0.999 \\
  Sample Splitting & 0.5 & 0.05 \\
  Pre-validation & 0.3 & 0.08 \\
  Selective Inference  & 0.5 & 0.05 \\
  Data Carving & 0.576 & 0.05 \\
\hline
\end{tabular}
\label{tab:null}
\end{table}

As we can see, this is a nontrivial problem because just performing the naive t-test (pretending that $\hat{y}$ is given) leads to a huge type I error.  This is unsurprising, because we used the lasso to find a $\hat{y}$ that is close to $y$.  One obvious method that could be used to address this problems is sample splitting -- using half of the data to generate our second stage hypothesis and using the other half to test it.  We discuss sample splitting along with a related method called pre-validation in Section~\ref{sec:prob}.  Pre-validation seeks to do a better job using the available information in the problem by using all of the data in both stages.  While this makes rejections more likely, pre-validation fails to protect type I error.

Selective Inference is a relatively new technique that gives exact p-values for many tests that are conducted on data after some sort of selection procedure.  In our case, the selection procedure is defined by the columns of $X_2$ that are selected by the fitter.  By conditioning on the selection, we are able to use all of our data in both stages and still control type I error.  In Section~\ref{sec:posi} we discuss how selective inference works.  Data carving is a type of selective inference when not all of the data is used in the first stage (in this way it is somewhat similar to sample splitting); reserving data can increase the power of the second stage test.  In Section~\ref{sec:siip} we will show how to use selective inference techniques (including data carving) in order to test the inclusion of $\hat{y}$.  

In Section~\ref{sec:iif} we show how to use selective inference to perform a test on the inclusion of $X_{\hat{E}}$.  While this problem may not be as directly tied to internal inference as the inclusion of $\hat{y}$, it is still potentially interesting.  Additionally, in this case it is possible to derive an analytical solution.

We apply our method to a real biostatistics problem in Section~\ref{sec:realdat} to show how it can be used in practice.  Finally, we conduct several simulations in Section~\ref{sec:sims} to verify that our method properly controls the level of the test and to evaluate our ability to find true signals.

\section{Existing Approaches}
\label{sec:prob}

In the case where $\hat{y}$ does not vary with $y$, then we could just use a classical z/t-test (in the case where the variance of $\epsilon$ is /isn't known) to test the inclusion of $\hat{y}$.  However, as implied by $\hat{y} \equiv f(X;X,y)$, $\hat{y}$ will typically vary with $y$.  This means that $y$ appears on both sides of our model equation, and the classical tests will no longer be valid.


\subsection{Sample Splitting}

When the sample size is large, there is a simple fix.  We can first split our dataset into two parts: a training set and a test set.  $\hat{y}$ can be learned on the training set, and then the
predictive power of the learned signature can be studied on the 
test set\citep{acox1975note}.  Since the training set and test set are presumably independent, we can use standard statistical approaches like a t-test for including the signature in a regression of $y$ on $Z$.  Formally, sample splitting proceeds as follows:

\begin{enumerate}
  \item Partition the dataset into two parts, $(X_1, Z_1, y_1)$ and $(X_2, Z_2, y_2)$ by putting $n_1$ of the observations into a training set and reserving the $n_2 = n - n_1$ remaining observations into a test set.
  \item Let $\hat{y_2} = f(X_2; X_1, y_1)$ for some fitter of interest $f$.
  \item Perform the t-test for the inclusion of $\hat{y_2}$ in a regression of $y_2$ on $Z_2$.
\end{enumerate}

Similarly, if we want to test the inclusion of $X_{\hat{E}}$, we would find the $\hat{E}$ implied by $(X_1, y_1)$, and then perform the F-test for including $(X_2)_{\hat{E}}$ in a regression of $y_2$ on $Z_2$.

This idea, referred to as sample splitting, is a commonly taught in classes on statistics as a way to conduct ``fair'' analyses.  Sample splitting is a common approach in biostatistics (where the withheld set of data is sometimes referred to as a validation set) as it provides a sort of internal reproduction of whatever result is implied by the training set \citep{abuhlmann2014high}.  The popularity of sample splitting is probably due to it's generality, ease of use, and validity.  That said, there are some disadvantages to sample splitting.  First, by only using the training set to fit a model for the internal predictor, we limit our ability to find a good internal predictor.  Second, because we only use the test set to test the significance of $\hat y$, it is possible that we are ignoring some leftover information from the training set.  \citet{afithian2014optimal} show that sample splitting is inadmissible in the situation we are considering (testing regression coefficients) for precisely this reason.

\subsection{Pre-Validation}
\label{sec:pre-val}

To address the issue of sample splitting not using all of the available data to find an interesting internal predictor, \citet{atibshirani2002pre} proposed the idea of {\em pre-validation} -- a  method that was already in 
use in the biomedical literature, but  without formal justification \citep{avan2002gene}.  Essentially, the idea is to split the observations into folds, as in $K$-fold cross validation.  For each fold, a model is learned on the rest of the data predicting the outcome from the biomarkers.  Then, a prediction is formed for the response on that fold.  The predictions for each of the folds are combined together to form the \emph{pre-validated} predictor for the outcome based on the  biomarkers.  Namely:

\begin{enumerate}
  \item Partition the dataset into $K$ folds -- $(X_1, Z_1, y_1)$, $(X_2, Z_2, y_2)$, \dots, $(X_K, Z_K, y_K)$.
  \item For each $i = 1, 2, \dots K$: let $\hat{y}_i = f(X_i; X_{[-i]}, y_{[-i]})$ where the negative indexing indicates removal of that fold.
  \item Test the inclusion of $\hat{y} = (\hat{y}_1, \hat{y}_2, \dots \hat{y}_K)$ in a regression of $y$ on $Z$.
\end{enumerate}

Pre-validation tries to arrange a fair comparison at the second stage, because each prediction is made independent of the response value for that observation.  At the same time, relative to sample splitting, more observations can be used in the formation of the pre-validated predictor.
\citet{ahofling2008study} studied the theoretical properties of
pre-validation.  They showed that the pre-validation approach leads to
invalid p-values and recommend using a permutation null to account for
this fact.  One challenge of the permutation null is the computational
cost with which it is associated (one must repeat an algorithm that
already involves a cross validation many times).  In practice, while
pre-validation is used in modern biology, the examples we examined used
the naive method instead of the permutation null
\citep{aspahn2010expression, asegura2010melanoma}.

Note that pre-validation does not provide a way to test the inclusion of $X_{\hat{E}}$.  It was designed for internal inference, and does not approach the related problem.

\section{Selective Inference}
\label{sec:posi}

Here we briefly summarize the theory for inference after model selection of \citet{afithian2014optimal}.  They create a theory for analyzing statistical hypotheses when those hypotheses are data driven.  This is often the case in general settings where statistics is used, and is certainly the case in our situation (where the hypothesis we test will depend on $\hat y$ or $X_{\hat{E}}$).  Essentially, instead of controlling the type I error, they create a framework for controlling the \emph{selective type I error rate}.  

If we want to test hypothesis $H_0$ under model $M$, then controlling type I error would mean to limit $P_{M, H_0}(\textrm{reject } H_0)$.  By contrast, controlling selective type I error would mean limiting $P_{M, H_0}(\textrm{reject } H_0 | (M, H_0) \textrm{ selected})$.  By conditioning on the selection of model and hypothesis, we prevent ourselves from using the same information to both generate and test a hypothesis.  In the case of sample splitting, the data used to select $(M, H_0)$ is independent of the data used to test $H_0$, so controlling type I error is equivalent to controlling selective type I error.

One simple example where the importance and utility of selective inference is apparent is in determining the bias of a coin.  Imagine if you flipped a coin 9 times and it came up heads 8 times.  You might hypothesize that the coin is biased towards heads and decide to conduct a test for the null hypothesis of a fair coin against an alternative that the coin is biased towards heads.  We can get a p-value for this test from the binomial distribution function.  The probability of a fair coin having 8 or more heads is around .01.  The test we performed controls type I error, so our procedure is valid from the view of classical statistic.  That said, readers should be disappointed by the use of the one sided test here.  After all, we only decided to test that the coin was biased towards heads after we saw our sample.  If we had seen 8 tails we probably would have ended up testing a bias towards tails.  One way to resolve the issue in our biased coin example is to use the two sided test.  

Another way to resolve the issue above would be by applying selective inference.  In this case, we have data $X$ equal to the number of heads observed in 9 flips.  If $X \ge 5$, then we will test the null hypothesis ($H_0$) that the probability of heads, $p$, equals 1/2 in the model ($M$) $p \in [1/2, 1]$.  If $X \le 5$ we will test the same null hypothesis for the model $p \in [0, 1/2]$.  So, instead of using the binomial distribution function to look at $P_{p=.5}(X \ge 8)$, we should test conditional on us choosing the model $p \in [1/2, 1]$ (equivalent in this case to choosing which tail to test).  This means our p-value will be $P_{p=1/2}(X \ge 8 | X \ge 5)$.  This gives a value around .02 and it is easy to verify that due to the symmetry of the binomial distribution with $p=1/2$ the selective test described above is equivalent to the two sided test.

In this way, selective inference formalizes the issue that people have with the use of a one sided test above.  The issue is not that the test fails to control type I error (it does control type I error).  The issue is that it fails to control type I error conditional on the test actually being run.

\citet{afithian2014optimal} demonstrate that collections of selective tests will achieve a long run control over the frequentist error rate: $\frac{\# \textrm{false rejections}}{\# \textrm{true nulls selected}}$.  They refer to this property as Discipline-Wide Error Control and it serves as a response to the common critique of type I error control that aggregating only significant results into journals will lead to a frequentist error rate far above whatever level at which the type I error is controlled.  They further note that no such result exists for false discovery rate or family-wise error rate, making selective inference an attractive alternative viewpoint to the way multiple testing is typically treated.

\section{Selective Inference on Internal Predictors}
\label{sec:siip}

Here, we lay out our theory for a selective test on the second stage hypothesis:

\[
  H_{0,M}:\theta_M = \theta_0,\;\;\;y = \theta_M\hat{y}   + Z\gamma + \epsilon
\]
where $\epsilon$ is an iid normally distributed noise vector.  We have added $M$ subscripts to demonstrate that the null hypothesis here is dependent on the selection procedure.

As discussed before, we focus on cases where the selection implies some sort of variable selection.  Namely, on the event $\hat{M} = M$ there is a function $f_M$ such that $\hat{y} = f(X; X, y) = f_M(X_E; X, y)$ where $E$ is a subset of the columns of $X$ implied by $M$.

We restrict our attention to \emph{polyhedral selection procedures} -- ones where the event $\hat{M} = M$ is equivalent to $A_My \le b_M$ for some appropriately chosen $A_M$ and $b_M$.  This includes many interesting cases including marginal screening, forward stepwise regression, and the lasso to name a few \citep{alee2014exact,aloftus2014significance,alee2013exact}.

\subsection{Internal Predictor Depends on $y$ Only Through Selection}
\label{sec:notar}

In the case where our internal predictor is independent of $y$ conditional on $\hat{M} = M$ then our problem has already been addressed in the selective inference literature.  In this case, we are just testing the inclusion of $f'_{M}(X_E; X) \equiv f_M(X_E; X, y)$ (a constant conditional on $\hat{M} = M$) in a regression of $y$ on $Z$.  Thus, we want to test if the first element of $(\tilde{Z}^T\tilde{Z})^{-1}\tilde{Z}^Ty$ is 0, where $\tilde{Z} = (\hat{y} \; Z)$.  \citet{ataylor2014post} and \citet{alee2013exact} show that if $y$ is normally distributed, this statistic conditional on $M$ and the sufficient statistics for the nuisance parameters ($Z^Ty$) will have a truncated normal distribution with truncation points that can be found from $A_M$, $b_M$, and $\text{var}(y)$.  Thus, we can use the truncated normal distribution to derive p-values for testing the inclusion of $\hat{y}$'s of this form.

While this may seem like a simple case, it encompasses several interesting approaches for building an internal predictor:
\\

{\em Example 1: Predictor most correlated with $y$.}
One simple selection procedure uses just the predictor with highest absolute
correlation with $y$ (WLOG assume it is $\xs$), so that  $\xs$ is our internal predictor
of $y$.  Thus, in this case  the selection procedure is equivalent to
$|\xs^Ty| \ge |x_i^Ty| \;\forall i > 1$.  We can rewrite this as $(\textrm{sign}(\xs^Ty)\xs^T - \textrm{sign}(x_i^Ty)x_i^T)y \le 0 \;\forall i$.  This has the form $Ay\leq b$ where $b = 0$ and $A$ has rows $(\textrm{sign}(\xs^Ty)\xs^T - \textrm{sign}(x_i^Ty)x_i^T)$ for $i = 2, 3, \dots p$.  Note that our selection event here, $M$, is just the signs of $x_i^Ty$ (if we hadn't sorted the columns of $X$ it would also include the index with the highest correlation).  Since $\hat{y} = \xs$ which is constant given $M$, we can use the method described above to perform valid inference on including $\hat{y}$ in a regression of $y$ on $Z$.
\\

{\em Example 2: Fixed linear combinations after marginal screening.}
In addition to the example where we select just one predictor, it is also possible to explain the selection of the top $k$ predictors of $y$ as an affine selection event (WLOG assume they are the first $k$ predictors).  In this case, we have $|x_j^Ty| \ge |x_i^Ty| \;\forall i > k, j \le k$.  Again, our selection event is the signs of $x_i^Ty$.  Then, as long as we have an $\hat{y}$ that is constant with respect to $y$ given $M$, we can use the above tools.  For example, maybe we take the average or first principal component of the top $k$ correlated predictors and use that as an internal predictor for which we want a p-value.
\\

That said, there are also interesting cases where $\hat{y}$ varies with $y$ conditional on $M$.  For example, consider building an internal predictor by fitting a Lasso regression of $y$ on $X$ with
a fixed regularization parameter $\lambda$, giving sparse coefficient vector $\hat{\eta}$.  Letting $E = \{i: \hat{\eta}_i \ne 0\}$, and $z = \textrm{sign}(\hat\eta)$, we can
express this as a selection event where $M = \{E, z_E\}$ (namely, the
selected variables and the signs of their coefficients).  The Lasso
was a motivating example for \citet{alee2013exact}, who give the formulation for $A_M$ and $b_M$, but a problem remains that prevents us from using their method.  The $\hat y$ fit from a Lasso is not constant with respect to $y$ given $M$ (which implies the $\eta$ is not constant); It is equal to $\hat y = L_My + \zeta_M$ where $L_M$ is $X_E(X_E^TX_E)^{-1}X_E^T$ and $\zeta_M$ is $-\lambda (X_E^T)^{\dagger}z_E$.  This means that while we could form valid p-values for including any of the variables in $X_E$ to a regression of $y$, we could not use the above method to form a p-value for the inclusion of $\hat y$.

\subsection{Internal Predictor an Affine Function of $y$ Conditional on Selection}
\label{sec:intt}

We now present a method for assessing the significance of an internal predictor when that predictor is an affine function of $y$.  This will extend the case where the internal predictor is constant conditional on the selected model.  While this may seem complicated, we note that it is only slightly different from deriving the t-test for ordinary regression.  This is because our $\hat{y}$ is not actually a function of the parameter of interest.  The modification comes in adjusting for our selection event, and that is handled by conditioning.

On the event $\hat{M}=M$ we suppose we have a internal predictor
\be
\hat{y} = f_M(X_E; X, y) = L_My + \zeta_M
\ee
and our statistical model is parameterized by $(\theta_M,\gamma,\sigma^2)$:
$$
y  = \theta_M \hat{y} + Z \gamma + \epsilon
$$
where $\epsilon \sim N(0, \sigma^2 I)$.

Note that this model is the same as
\be
\begin{aligned}
y = (I - \theta_ML_M)^{-1} (\theta_M\zeta_M + Z \gamma + \epsilon ) \sim N \big(&(I - \theta_ML_M)^{-1} (\theta_M \zeta_M + Z \gamma),\\
& \sigma^2 (I - \theta_ML_M)^{-1}(I - \theta_ML_M)^{-T}\big)
\end{aligned}
\ee
(as long as $\theta_M$ not a reciprocal eigenvalue of $L_M$).

We want to test $H_{0,M}:\theta_M=\theta_0$. Under this hypothesis, unconditionally, we have
\be
R_1(y,\theta_0) \equiv (I - P_Z) (y - \theta_0\hat{y}_M) = (I - P_Z)\epsilon \sim N(0, \sigma^2 (I - P_Z))
\ee
where $P_Z = ZZ^{\dagger}$. Note that since $R_1$ has been defined to be the part of $\epsilon$ not explained by $Z$, it's distribution does  depend on our nuisance parameter $\gamma$. 

We can perform a test using the score statistic for this distribution conditional on $P_Zy$ and $\hat{M}=M$.  This statistic ends up being 

\be
T_t(y) = \frac{R_1(y, \theta_0)^T(I-P_Z)\hat{y}}{\sqrt{\hat{y}^T(I-P_Z)\hat{y}}}.
\ee
Note that this is the same statistic that would be used in the t-test for testing the inclusion of $\hat{y}$ in a regression of $y$ on $Z$ if $\hat{y}$ was an external predictor.  The derivation is identical.

As a reminder, the selection event is affine, i.e.
\be
\left\{y: \hat{M}(y)=M\right\} = \left\{y : A_My \leq b_M\right\}.
\ee
We then break $y$ into two independent parts, one that independent of $R_1$ and one that is explained by $R_1$:

We write
\be
\begin{aligned}
y &=& y - E_0(y | R_1(y,\theta_0)) &\;\;\;+& &E_0 (y | R_1(y,\theta_0))\\
&&\textrm{part independent of } R_1 &&& \textrm{part explained by } R_1.
\end{aligned}
\ee
where
\be
\begin{aligned}
E_0(y | R_1) &= \textrm{Cov}\left(y, R_1\right)\textrm{Var}\left(R_1\right)^{-1}R_1 \\
&= (I - \theta_0L_M)^{-1}R_1.
\end{aligned}
\ee

By construction, the two pieces of $y$ above are independent and we can rewrite our selection event as
\be
\left\{y: A_M \left((I - \theta_0L_M)^{-1}R_1(y,\theta_0) + \Delta(y,\theta_0) \right) \leq b_M \right\}
\ee
where $\Delta(y,\theta_0) = y - (I - \theta_0L_M)^{-1} R_1(y,\theta_0)$ is independent of $R_1(y,\theta_0)$.

We now condition on $\{\Delta(y,\theta_0)=\delta\}$ and the selection event is expressed as
\be
\left\{y: A_M (I - \theta_0L_M)^{-1}R_1(y,\theta_0) + A_M \delta \leq b_M \right\}
\ee

If we know $\sigma^2$ we can sample from this selection event and, say, use our test statistic $T(R_1(y,\theta_0))$ to generate
 a reference distribution to test $H_0$. Here, ``sample from this
 selection event'' means to draw $R_1 \sim N(0, \sigma^2 (I-P_Z))$
keeping only those that fall into the selection event.  This sampling could be done using an accept reject scheme, or using other samplers such as hit and run \citep{agolchi2014sequentially}.  In practice, for sufficiently large problems the accept reject sampling is infeasible, so we use a hit and run sampler to draw samples.

If we do not know $\sigma^2$, then we might also condition on $\|R_1(y,\theta_0)\|_2=l$ which leaves the only remaining
randomness in the unit vector $R_1(y,\theta_0)/\|R_1(y,\theta_0)\|_2$. The selection event is now
\be
\left\{y: \|R_1(y,\theta_0)\|_2=l, A_M (I - \theta_0L_M)^{-1} R_1(y,\theta_0) + A_M \delta \leq b_M \right\},
\ee
which can be similarly sampled. 

\subsection{P-values for General T-test}

Additionally, if we are only interested in testing $\theta_M = 0$, we can perform a test using a general fitter for $y$: $\hat{y} = f(X;X,y)$. This is possible because under this null hypothesis, unconditionally:
\[
\begin{aligned}
y = Z \gamma + \epsilon  \sim N \big(Z \gamma, \sigma^2 I\big).
\end{aligned}
\]

Thus, $E_0(y | R_1) = R_1$ and we can do sampling as before, but now based on the selection event:

$$
\left\{y: A_MR_1(y,\theta_0) + A_M \delta \leq b_M \right\}.
$$

\subsection{Data Carving}

Data Carving is a term used to describe selective inference methods that only use a subset of the data to perform selection.  We can do a data carving version of our method as follows:

\begin{enumerate}
  \item Partition the dataset into two parts -- $(X_1, Z_1, y_1)$ and $(X_2, Z_2, y_2)$ -- by putting $n_1$ observations into one part and retaining the other $n_2$ observations for the second part.
  \item Let $\hat{y} = f(X; X_1, y_1)$ for some fitter of interest $f$.
  \item Calculate the $A_M$ and $b_M$ that represent the selection event implied by $(X_1, y_1)$.
  \item Add $n_2$ columns of 0 to $A_M$.
  \item Perform the selective inference method described above on the full dataset (concatenate the two parts back together) with $\hat{y}$, $b_M$, and our modified $A_M$.
\end{enumerate}

By making the slight change to $A_M$ described above, we ensure that only the first part of $y$ influences the selection event.  Data carving resembles sample splitting but allows us to use the full dataset in the second stage.  By conditioning on the selection event we avoid reusing information from the first part, but the fact that we get to use the additional information contained in that part means that we should perform better than sample splitting.  In fact, \citet{afithian2014optimal} show that sample splitting is often inadmissible through the existence of data carving.

\section{Testing the inclusion of $X_{\hat{E}}$}
\label{sec:iif}

Above, we present a method for calculating a selective t-statistic for the inclusion of $\hat{y}$ in a regression of $y$ on $Z$.  This makes sense for the typical internal inference problem.  For some post selection selection settings though, we might want to do a related test of the inclusion of the columns of $X_{\hat{E}}$ in a regression of $y$ on $Z$.  

On the event $\hat{M} = M$, assume we have the following statistical model parameterized by $(\beta_M,\gamma,\sigma^2)$:

\be
y = X_E\beta_M + Z \gamma + \epsilon
\ee
where $\epsilon \sim N(0, \sigma^2 I)$.  

We want to test $H_{0,M}:\beta_M=\beta_0$. Under this hypothesis, unconditionally, we have
\be
R_1(y,\beta_0) = (I - P_Z) (y - X_E\beta_0) = (I - P_Z)\epsilon \sim N(0, \sigma^2 (I - P_Z)).
\ee

Define also
\be
R_2(y,\beta_0) = (I - P_M)(y) = (I - P_M) \epsilon \sim N(0, \sigma^2 (I - P_M))
\ee
where the distribution holds unconditionally and $P_M$ is projection onto the column space of 
$$
\tilde{X}_M = \begin{pmatrix}
X_E & Z 
\end{pmatrix}.
$$
Note that $R_2(y,\beta_0)=(I-P_M)R_1(y,\beta_0)$.

Therefore, unconditionally
\be
T_F(R_1(y,\beta_0)) = \frac{\frac{\|R_1(y,\beta_0)\|^2_2 - \|R_2(y,\beta_0)\|^2_2}{\text{Tr}(P_M-P_Z)}}{\frac{\|R_2(y,\beta_0)\|^2_2}{\text{Tr}(I-P_M)}} \sim F_{\text{Tr}(P_M-P_Z), \text{Tr}(I-P_M)}.
\ee
We will construct a test based on the selective distribution of $T_F$.  Note that this is just the derivation of an F-test for including variables in ordinary regression.

In the simplest scenario, the selection event is affine, i.e.
$$
\left\{y: \hat{M}(y)=M\right\} = \left\{y : A_My \leq b_M\right\}.
$$
We then break $y$ into two independent parts, one that is explained by $R_1$ and one that is independent of $R_1$ (and thus $T_F$):

We write
\be
\begin{aligned}
y &=& y - E_0(y | R_1(y,\beta_0)) &\;\;\;+& &E_0 (y | R_1(y,\beta_0))\\
&&\textrm{part explained by } Z &&& \textrm{part not explained by } Z.
\end{aligned}
\ee
where
\be
\begin{aligned}
E_0(y | R_1) &= \textrm{Cov}\left(y, R_1\right)\textrm{Var}\left(R_1\right)^{-1}R_1 \\
&= R_1.
\end{aligned}
\ee

By construction, the two pieces of $y$ above are independent and we can rewrite our selection event as
\be
\left\{y: A_M \left(R_1(y,\beta_0) + \Delta(y,\beta_0) \right) \leq b_M \right\}
\ee
where $\Delta(y,\beta_0) = y - R_1(y,\beta_0)$ is independent of $R_1(y,\beta_0)$.

We now condition on $\{\Delta(y,\beta_0)=\delta\}$ and the selection event is expressed as
\be
\left\{y: A_M R_1(y,\beta_0) + A_M \delta \leq b_M \right\}
\ee

If we know $\sigma^2$ we can sample from this selection event and, say, use our test statistic $T(R_1(y,\beta_0))$ to generate
 a reference distribution to test $H_0$. Here, ``sample from this
 selection event'' means to draw $R_1 \sim N(0, \sigma^2 (I-P_Z))$
keeping only those that fall into the selection event (of course this does not need to be accept reject -- could be hit and run or other
sampler).  

If we do not know $\sigma^2$, then we might also condition on $\|R_1(y,\beta_0)\|_2=l$ which leaves the only remaining
randomness in the unit vector $R_1(y,\beta_0)/\|R_1(y,\beta_0)\|_2$. The selection event is now
\be
\left\{y: \|R_1(y,\beta_0)\|_2=l, A_M R_1(y,\beta_0) + A_M \delta \leq b_M \right\},
\ee
which can be similarly sampled. 

\subsection{A non-sampling test for $\beta_0=0$}
\label{sec:nonsamp}

In this case, we have the ``null model'':
\be
y = Z\gamma + \epsilon, \qquad \epsilon \sim N(0, \sigma^2 I).
\ee

The sampling distribution we draw from is $y \sim N(0, \sigma^2I)$ subject to the polyhedral constraints
\be
\left\{A_M (I-P_Z)y  + A_M \delta \leq b_M \right\}
\ee
with $\delta$ the observed value of $P_Zy$.

Suppose we take the last test which conditions on $\|R_1(y,\beta_0)\|_2$ and $P_Zy$ and also condition on the unit vectors
\be
V_N = \frac{(P_M-P_Z)y}{\|(P_M-P_Z)y\|_2}, \qquad V_D = \frac{(I-P_M)y}{\|(I-P_M)y\|_2}.
\ee
Then $R_1(y,\beta_0)/l=\frac{(cT_F)^{1/2}}{(1 + cT_F)^{1/2}}\cdot V_N + \frac{1}{(1 + cT_F)^{1/2}} \cdot V_D$ where $c = \frac{\text{Tr}(P_M-P_Z)}{\text{Tr}(I-P_M)}$
and the only remaining variation is the statistic $T$ above. The selection event can be written as
\be
\left\{T_F \geq 0: A_M \left(\frac{(cT_F)^{1/2}}{(1 + cT_F)^{1/2}}\cdot V_N + \frac{1}{(1 + cT)^{1/2}} \cdot V_D\right) \cdot l + A_M \delta \leq b_M \right\}.
\ee
with $\delta$ the observed value of $P_Zy$.

We can then use the appropriate $F$ distribution truncated to the above set as a reference distribution.  A proof of this claim and details about how to find the truncation set for $T_F$ can be found in Appendix~\ref{app:nonsamp}.




\section{Real Data Example}
\label{sec:realdat}

In order to demonstrate the use of our method for internal inference, we perform a similar analysis to the one that \citet{atibshirani2002pre} and \citet{ahofling2008study} performed on the dataset from \citet{avan2002gene}.  This dataset consists of 78 patients with breast cancer that have been split into 44 patients with a good prognosis, and 34 patients with a poor prognosis.  For each patient there are 4,918 gene expression measurements as well as the following clinical predictors: tumor grade, estrogen receptor status, progestron receptor status, tumor size, patient age, and angioinvasion.  The goal of the analysis is to determine whether a gene signature learned on the gene expression data will add value beyond that of the clinical predictors in predicting prognosis.

In those analyses, the constructed predictor is based on a complicated method that involves screening variables against the response and then measuring correlations with the class centroids to build a decision rule.  This highlights one relative advantage of the pre-validation and sample-splitting approaches; they can be used on arbitrary methods of building a constructed predictors.  In order to make a comparison between methods, we instead explore building an internal predictor using a Lasso regression with our response coded as $y_1 = 1$ for patients with good prognosis and $y_i = -1$ for those with a poor prognosis.  Note that the work of \citet{atian2015asymptotics} suggests that selective inference techniques are also possible for $\ell_1$ penalized logistic regression, but they rely on a normal approximation and thus are not exact.

We calculated p-values for the statistical significance for adding the
gene signature learned on the data ($\hat{y}$) to $Z$ in a model for $y$ using
several methods.  In order to select $\lambda$ we looked at the number of nonzero coefficients in the first stage model for several multiples of $2\sqrt{2\log p}$ and picked one that gave around 10-15 nonzero coefficients \citep{acandes2009near}.  Two of the methods we examine use all of the data to form a null hypothesis: our test and prevalidation.  The other two methods only use half of the data: sample splitting and a data carving version of our method.  Thus, it is important to note that while all of the p-values
are being reported in Figure \ref{fig:veerps} together, they are not necessarily directly
comparable.  This is because all of the tests may have differing null
hypothesis (depending on the selection event).  We also report the naive p-value from just assuming that $\hat{y}$ was an external predictor.  Although sample splitting does very well in this instance, the results for it and pre-validation are largely seed dependent.  It is difficult to read too much into the result of one analysis, and this is mainly intended to show how our selective inference approach could be applied.
\begin{figure}
\centering
\caption{Results of Analysis on Veer Data}
\begin{tabular}{ | l | l |}
\hline
  Method & p-value \\
\hline \hline
  Selective Inference & .127 \\
  Pre-validation  & 0.74 \\
  Naive T-test & 1.69e-07 \\
  Sample Splitting & 0.004 \\
  Data Carving  & 0.128 \\
\hline
\end{tabular}
\label{fig:veerps}
\end{figure}

\section{Simulations}
\label{sec:sims}

We also conducted a large simulation study in order to evaluate the
level and probability of rejection of selective inference and sample splitting under a setting where the truth is known.  In these simulations, we:

\begin{enumerate}
  \item Generate $X \in \mR^{n\times p_x}$ a matrix of independent standard normal random variables.  Center and scale $X$ to have unit variance.
  \item Generate $Z \in \mR^{n\times p_z}$ a matrix of independent standard normal random variables.  Center and scale $Z$ to have unit variance.
  \item Generate $y = X\beta + \frac{b_z}{\sqrt{pz}}Z\vec{1} + \epsilon$. Where $\beta$ is a vector with $p_{\textrm{real}}$ values of $\frac{b_x}{p_{\textrm{real}}}$ and the rest 0, $b_x$ and $b_z$ are parameters that control the relative importance of $X$ and $Z$, and $\epsilon$ is a noise vector of $n$ independent standard normal random variables.
\end{enumerate}

We chose to sample from the above model and not the model described in Section~\ref{sec:intt} mainly due to the fact that it would have given our method an unfair advantage in tests where we are looking at our ability to reject the selected null hypothesis.  This is due to the fact that it guarantees $y$ will be in the selection event, which is equivalent to saying that our first stage lasso will find the correct model.  Note that the model we do use to generate data has an equivalent null distribution.

\begin{figure}
    \centering
    \subfloat[Tests for inclusion of $\hat{y}$]{\includegraphics[width=.49\textwidth]{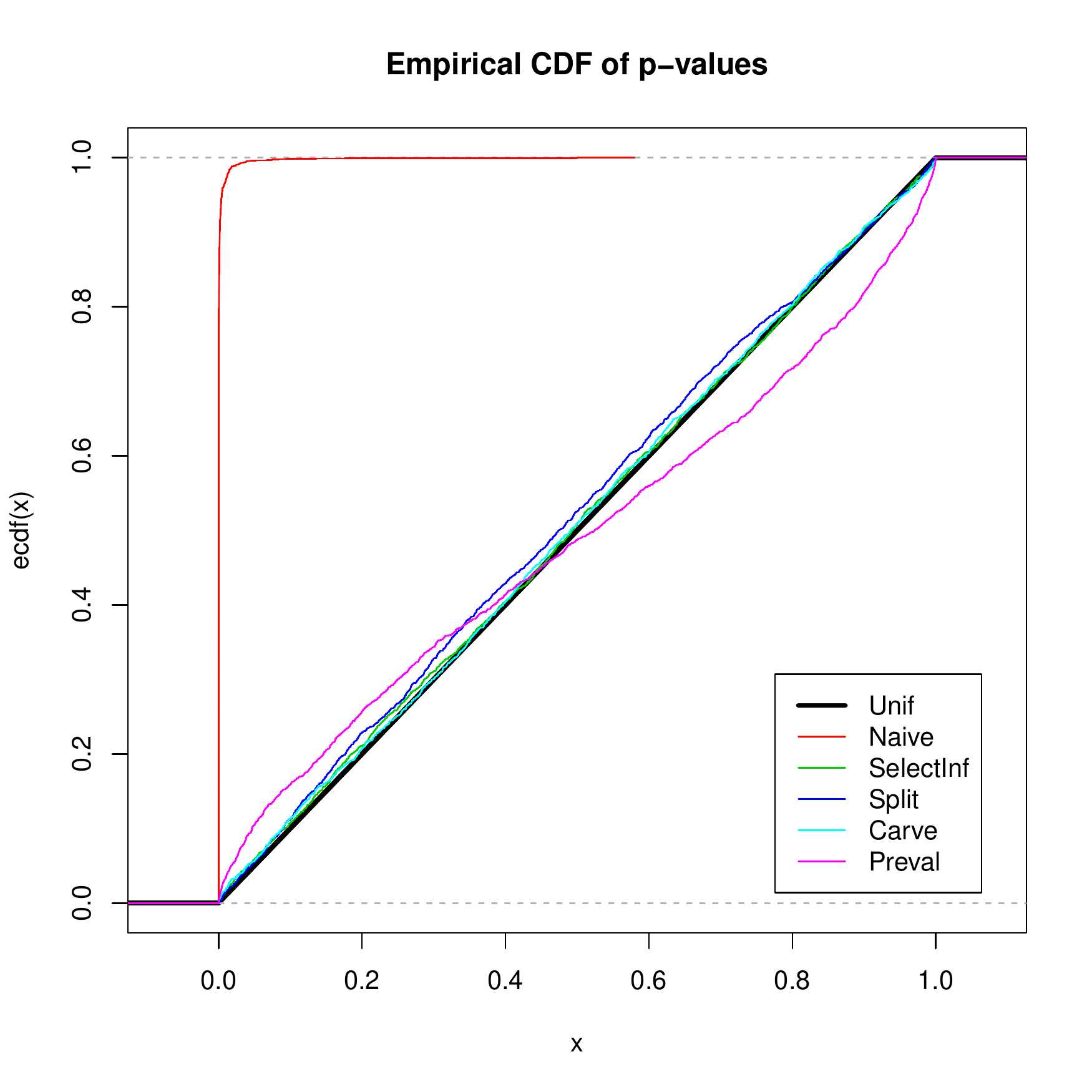}}
    \subfloat[Tests for inclusion of $X_{\hat{E}}$]{\includegraphics[width=.49\textwidth]{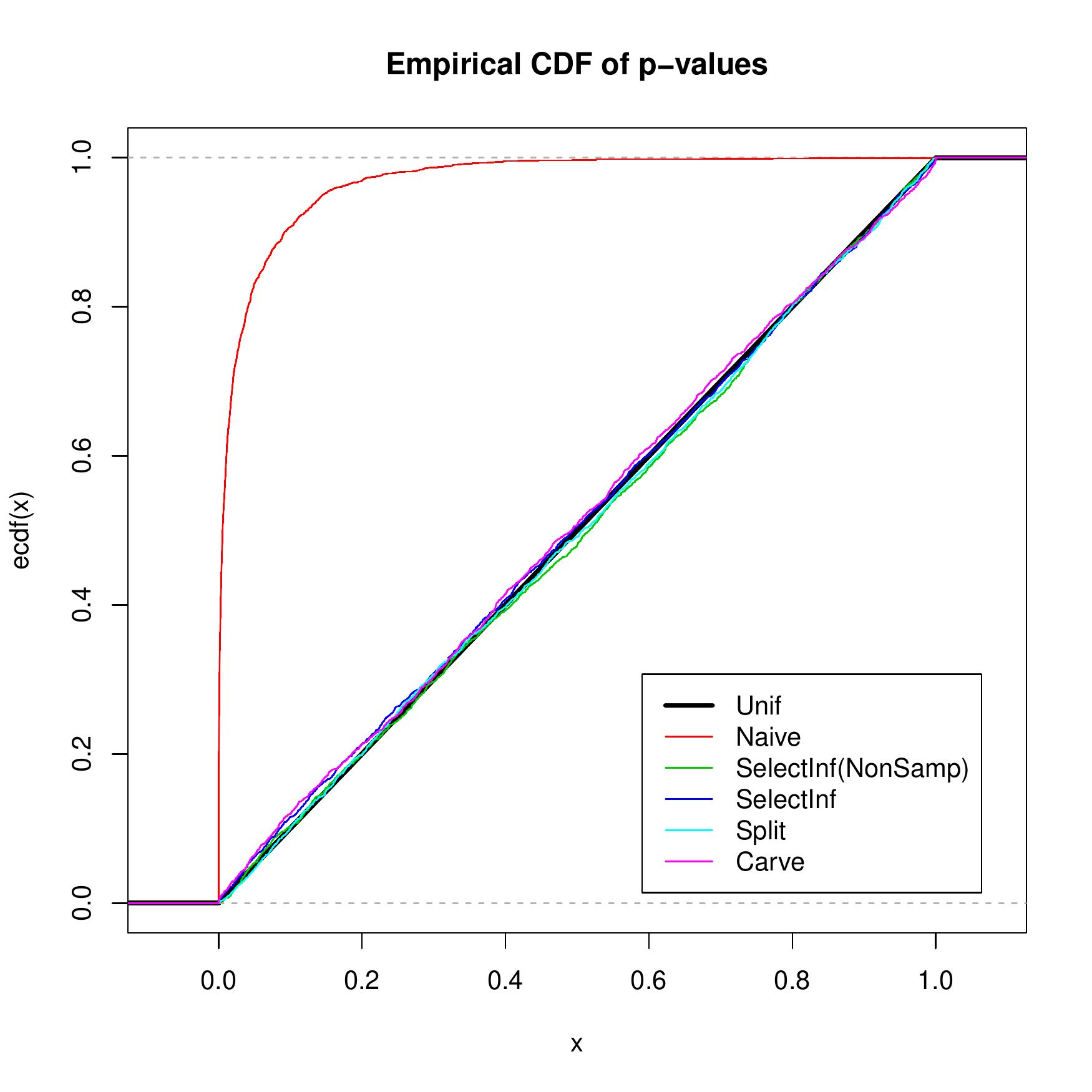}}
    \caption{Null Simulation}
    \label{fig:nullsim}
\end{figure}

Our first simulation uses $b_x = 0$ to test whether our method
properly protects selective type-I error.  We use $n = 50$, $p = 100$, $p_z =
5$, $p_{\textrm{real}} = 5$, $b_z = 1$ (giving
$\frac{b_z}{\sqrt{pz}}Z\vec{1}$ and $\epsilon$ the same variance).  We compare a naive t-test (as if $\hat{y}$ was given), our
method, sample splitting, pre-validation, and a data carving version of our test.  The data carving method uses the same splits as the sample splitting method.  Thus, some of our methods use all of the data to form the first stage model, and some only use half.  Both cases require a choice of lambda, and we selected a fixed value that gave approximately 10 nonzero coefficients in $\hat{\eta}$ for each model separately.  We also used the same simulations to test some of our methods that can be used on the test for the inclusion of $X_{\hat{E}}$: naive (as if $\hat{E}$ was given), our method, our non-sampling method, sample splitting, and a data carving version of our method.  Figure \ref{fig:nullsim} shows the p-values for 2000 runs of the above simulation.  As expected, all of the methods studied except for the naive and pre-validation methods properly protect type-I error.

\begin{figure}
  \centering
    \subfloat[Tests for inclusion of $\hat{y}$]{\includegraphics[width=.49\textwidth]{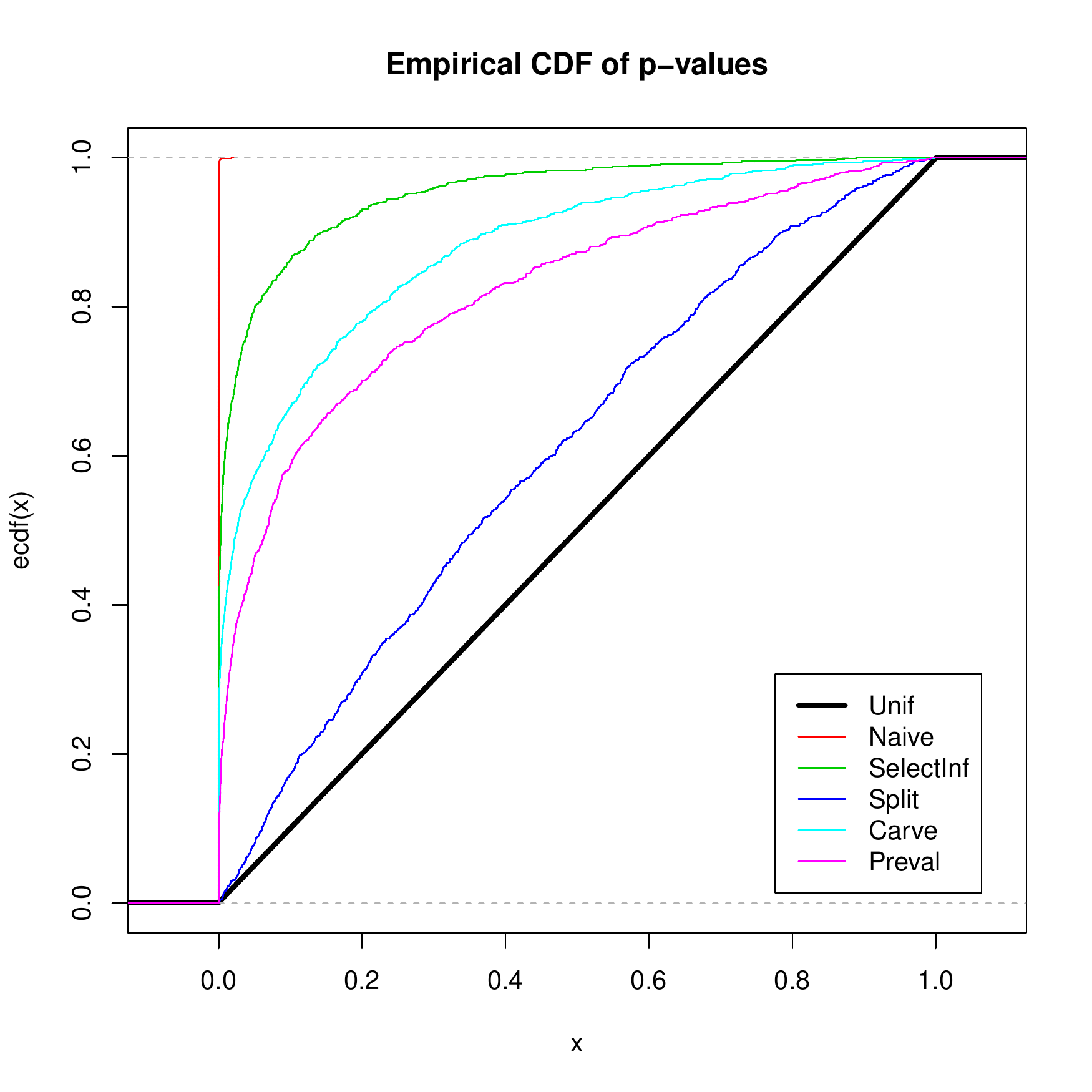}}
    \subfloat[Tests for inclusion of $X_{\hat{E}}$]{\includegraphics[width=.49\textwidth]{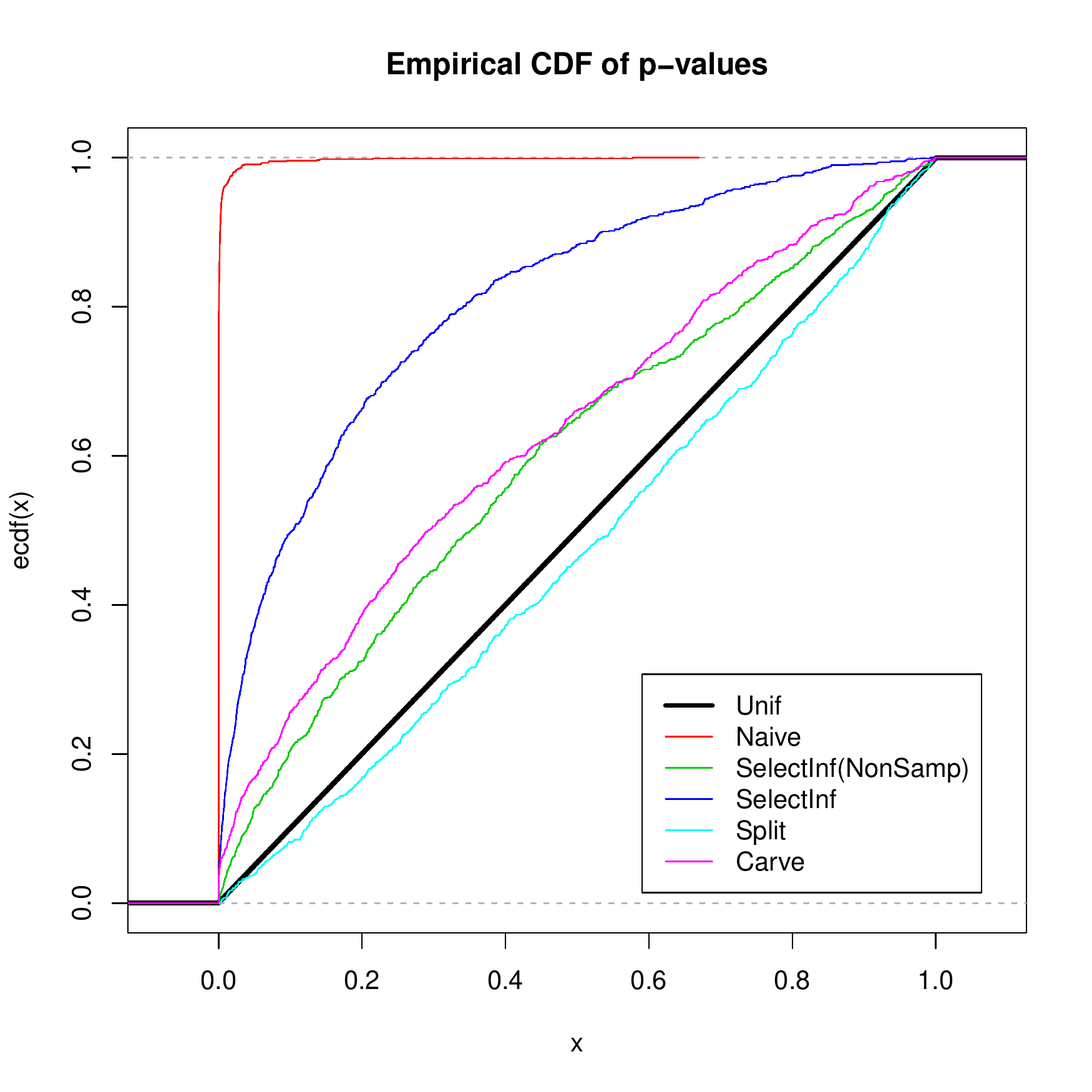}}
  \caption{Simulation with real effect.}
  \label{fig:powsim}
\end{figure}

We also ran the same simulation with $b_x = \sqrt{2}$ ($b_xX\beta$ explains half of the variance in $y$).  This allows us to see how well each of the methods identify true effects.  Again, $\lambda$ was selected separately for each method to ensure about 10 nonzero coefficients in the internal predictor.  The p-values from this simulation can be viewed in Figure \ref{fig:powsim}.  The methods that use the full dataset to form a hypothesis do a much better job of finding the columns of $X$ that have a true effect; while the methods that use the full dataset found 2.76 real coefficients on average, the methods that use half found 1.26.  We expect that a simulation where both methods find similar numbers of true coefficients would result in data carving outperforming the selective inference techniques.  Although the non-sampling version of our test for the inclusion of $X_{\hat{E}}$ is underpowered, it may still be convenient in cases where the sampling is not feasible.

\begin{figure}[h!]
  \centering
    \includegraphics[width=0.49\textwidth]{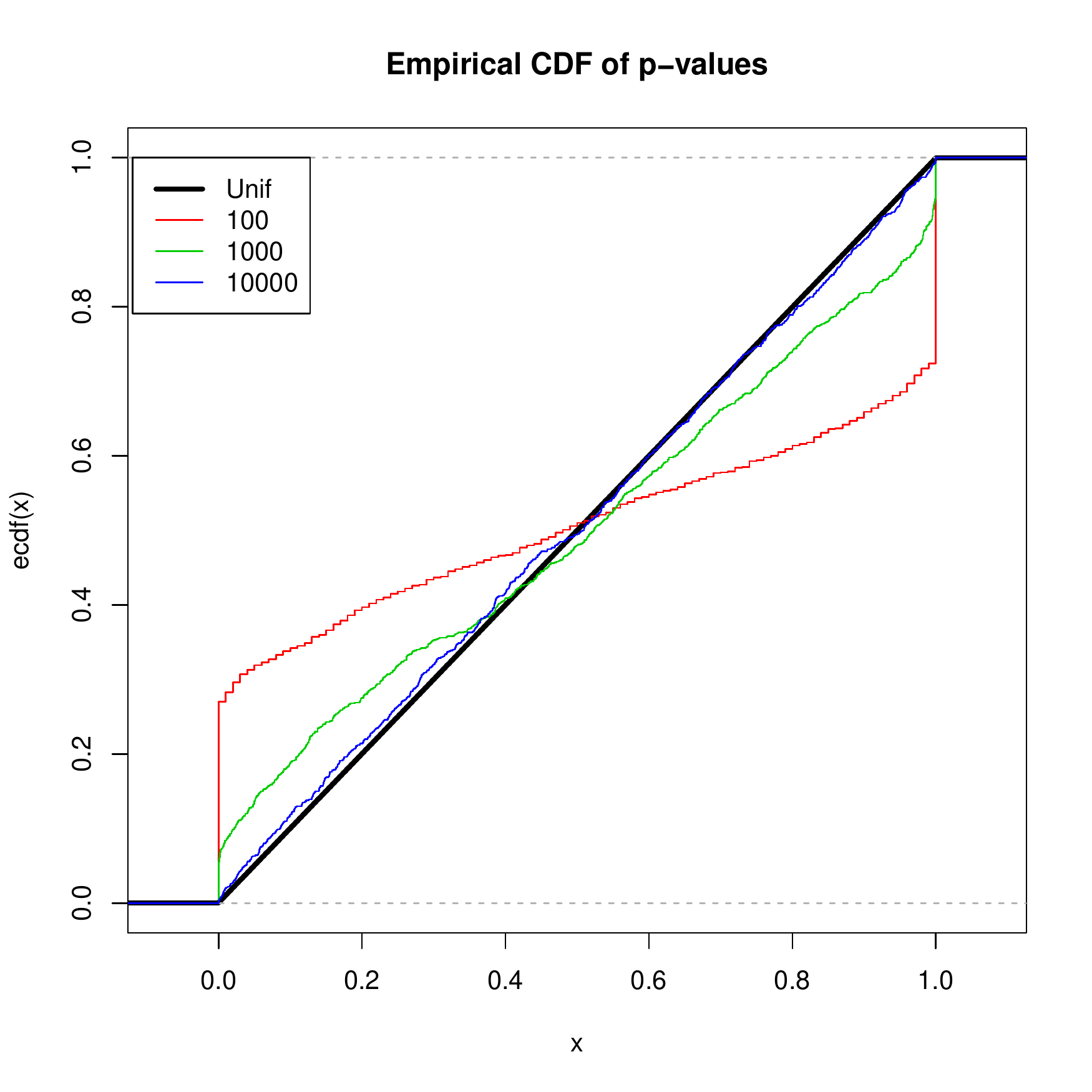}
  \caption{How many samples do we need to estimate p-value?}
  \label{fig:nsamps}
\end{figure}

Finally, we conducted a simulation to study the effect of not sampling enough when calculating the p-value of our internal t-test.  This is summarized in Figure \ref{fig:nsamps}.  The samples are taken from the selective inference method on the null sim for the inclusion of $\hat{y}$.  The different ecdf curves represent using different numbers of samples to compute each p-value.  Here, we see a serious issue that can arise in the sampling methods.  All sampling methods rely on having a sufficiently large sample to be valid.  If they are undersampled, then the p-values will have too much variance and will no longer be valid.  This can be especially challenging in weighted or autocorrelated samplers.

\section{Discussion}

In this paper we have explained how selective inference techniques can be used to
perform internal inference on high dimensional datasets.  We also discussed a similar problem of testing the inclusion of a set of predictors that have been found using a selection procedure.  This related problem has an analytical solution, which makes it convenient in cases where the sampling algorithms are not feasible.  Using selective inference 
has several advantages over competing approaches: it only requires one
pass through the modeling step (pre-validation requires a few for the
cross-validation and many more if using a permutation null), it
properly protects alpha, and it allows us to use more of the
information in the data to find an interesting hypothesis.  We also discussed how to implement data carving versions of our selective tools.  This means that in cases where it is easy to find good hypotheses, more data can be reserved for testing.

In addition to these advantages, there are several shortcomings in our
method that suggest areas for further research.  Most prominently, as
discussed in Section \ref{sec:realdat}, we have only implemented tools
for continuous responses.  Many interesting biological problems including our motivating example -- \citet{avan2002gene} -- have
responses that are binary
or survival data.  While it is possible to code binary data as
continuous (as we did), it would be more appropriate to develop tools
designed for that datatype.

Another major consideration when using our test in conjunction with the lasso is that our method is only correct for
fitting a lasso with fixed $\lambda$.  In practice, people typically
find interesting values of $\lambda$ using a cross validation or other
approaches.  Our method cannot currently adjust for this kind of
selection.

Despite the shortcomings, we feel our method is a useful tool for
modern biology.  These biomarker signatures are typically disregarded
until they can be successfully replicated, and our tool provides a way
of assessing which signatures are most likely to be successfully reproduced.

\appendix
\chapter{Details of Non-sampling Test}
\label{app:nonsamp}
Note that under $\beta_0=0$ we have $R_1 = (I - P_Z)y$, $R_2 = (I - P_M)y$, and $R_1-R_2$ is orthogonal to $R_2$ so $\|R_1 - R_2\|^2 = \|R_1\|^2 - \|R_2\|^2$,

\[
\begin{aligned}
& \frac{(cT)^{1/2}}{(1 + cT)^{1/2}}\cdot V_N + \frac{1}{(1 + cT)^{1/2}} \cdot V_D\\
& = \sqrt{\frac{\frac{\|R_1\|^2_2 - \|R_2\|^2_2}{\|R_2\|^2_2}}{\frac{\|R_1\|^2}{\|R_2\|^2}}} \cdot\frac{(P_M-P_Z)y}{\|(P_M-P_Z)y\|_2}  + \frac{1}{\sqrt{\frac{\|R_1\|^2}{\|R_2\|^2}}} \cdot \frac{(I-P_M)y}{\|(I-P_M)y\|_2} \\
& = \sqrt{\frac{\|R_1 - R_2\|^2_2}{\|R_1\|^2_2}}\cdot\frac{R_1 - R_2}{\|R_1 - R_2\|_2} + \sqrt{\frac{\|R_2\|^2}{\|R_1\|^2}}\cdot \frac{R_2}{\|R_2\|_2}\\
& = \frac{R_1 - R_2}{\|R_1\|} + \frac{R_2}{\|R_1\|} = U.
\end{aligned}
\]

We can rewrite our selection event as
\[
\left\{T \ge 0 :(cT)^{1/2}(A_MV_Nl) + A_MV_Dl + (1+cT)^{1/2}(A_M\delta - b) \le 0\right\}.
\]

Thus, we can interpret our selection event as limiting $cT$ to the intersection of sets of the form:
\[
\left\{x \ge 0 :g(x) \equiv q\sqrt{x} + r\sqrt{1+x} + s \le 0\right\}.
\]

We note that $\frac{\partial g}{\partial x} = \frac{q}{2\sqrt{x}} + \frac{r}{2\sqrt{1+x}}$ which means that $g$ has at most one local extremum.  This extremum can only occur if $\textrm{sign}(q) \ne \textrm{sign}(r)$ and $|r| > |q|$.  If it does exist, we can find all possible roots of $g$ by checking in the intervals $[0, \frac{q^2}{r^2 - q^2}]$ and $[\frac{q^2}{r^2 - q^2}, \infty]$, where $\infty$ can be replaced by a suitably large constant in practice.  If there is no local extremum, then we only need to look for roots on the interval $[0, \infty]$.  The precise set where $g(x) \le 0$ can then be found by looking at the signs of $g$ in the different intervals the roots separate.

\bibliographystyle{plainnat}
\bibliography{my.bib}{}

\end{document}